\documentclass{amsart}
\usepackage{amsmath,amsfonts,amssymb}
\usepackage{mathrsfs,geometry}
\usepackage[T1]{fontenc}
\usepackage[latin1]{inputenc}
\usepackage{url}

\newcommand{\iy}{\ensuremath{\infty}}
\newcommand{\R}{\ensuremath{\mathbf{R}}}
\newcommand{\C}{\ensuremath{\mathbf{C}}}

\newcommand{\N}{\ensuremath{\mathbf{N}}}

\newcommand{\E}{\ensuremath{\mathbf{E}}}

\newcommand{\<}{\ensuremath{\langle}}
\renewcommand{\>}{\ensuremath{\rangle}}
\renewcommand{\P}{\ensuremath{\mathbf{P}}}
\newcommand{\e}{\varepsilon}
\newcommand{\st}{\textnormal{ s.t. }}
\renewcommand{\leq}{\leqslant}
\renewcommand{\geq}{\geqslant}
\renewcommand{\phi}{\varphi}
\newcommand{\Id}{\mathrm{Id}}

\newcommand{\noi}{\noindent}

\newcommand{\mB}{\mathscr{B}}
\newcommand{\mU}{\mathcal{U}}
\newcommand{\mD}{\mathcal{D}}

\newcommand{\ket}[1]{| #1 \rangle}
\newcommand{\bra}[1]{\langle #1 |}
 
\DeclareMathOperator{\tr}{Tr}

\DeclareMathOperator{\conv}{conv}
\newcommand{\M}{\mathcal{M}}
\newcommand{\Ms}{\M_{\textnormal{sa}}}
\newcommand{\Mp}{\mathcal{M}_+}

\newtheorem{theo}{Theorem}
\newtheorem*{pr}{Proposition}
\newtheorem{lemma}{Lemma}

\newtheorem*{defn}{Definition}
\newtheorem{Atheo}{Theorem}
\newtheorem{Alemma}[Atheo]{Lemma}
\newtheorem*{Adefn}{Definition}

\begin{document}

\title{On almost randomizing channels with a short Kraus decomposition}
\author{Guillaume \sc{AUBRUN}}

\begin{abstract}
For large $d$, we study quantum channels on $\C^d$ obtained by selecting randomly $N$ independent Kraus operators according to a probability measure $\mu$ on the unitary group $\mU(d)$. When $\mu$ is the Haar measure, we show that for $N \succcurlyeq d/\e^2$, such a channel is $\e$-randomizing with high probability, which means that it maps every state within distance $\e/d$ (in operator norm) of the maximally mixed state. This slightly improves on a result by Hayden, Leung, Shor and Winter by optimizing their discretization argument. Moreover, for general $\mu$, we obtain a $\e$-randomizing channel provided $N \succcurlyeq d (\log d)^6/\e^2$. For $d=2^k$ ($k$ qubits), this includes Kraus operators obtained by tensoring $k$ random Pauli matrices. The proof uses recent results on empirical processes in Banach spaces.
\end{abstract}

\maketitle

\section{Introduction}

The completely randomizing quantum channel on $\C^d$ maps every state  to the maximally mixed state $\rho_*$. This channel is used to construct perfect encryption systems (see \cite{amtw} for formal definitions). However it is a complex object in the following sense: any Kraus decomposition must involve at least $d^2$ operators. It has been shown by Hayden, Leung, Shor and Winter \cite{hlsw} that this ``ideal'' channel can be efficiently emulated by lower-complexity channels, leading to approximate encryption systems. The key point is the existence of good approximations with much shorter Kraus decompositions. More precisely, say that a quantum channel $\Phi$ on $\C^d$ is $\e$-randomizing if for any state $\rho$, $\| \Phi(\rho)-\rho_* \|_{\iy} \leq \e/d$. The existence of $\e$-randomizing channels with $o(d^2)$ Kraus operators has several other implications \cite{hlsw}, such as counterexamples to multiplicativity conjectures \cite{winter}. 

It has been proved in \cite{hlsw} that if $(U_i)$ denote independent random matrices Haar-distributed on the unitary group $\mU(d)$, then the quantum channel
\begin{equation} \label{channel}  \Phi : \rho \mapsto \frac{1}{N} \sum_{j=1}^N U_i\rho U_i^\dagger \end{equation}
is $\e$-randomizing with high probability provided $N \geq  C d \log d/\e^2$ for some constant $C$. The proof uses a discretization argument and the fact that the Haar measure satisfies subgaussian estimates. We show a simple trick that allows to drop a $\log d$ factor: $\Phi$ is $\e$-randomizing when $N \geq C d / \e^2$, this is our theorem \ref{th1}.

The Haar measure is a nice object from the theoretical point of view, but is often too complicated to implement for concrete situations. Let us say that a measure $\mu$ on $\mU(d)$ is isotropic when $\int U\rho U^\dagger d\mu(U) = \rho_*$ for any state $\rho$. When $d=2^k$, an example of isotropic measure is given by assigning equal masses at $k$-wise tensor products of Pauli operators. 

The following question was asked in \cite{hlsw}: is the quantum channel $\Phi$ defined as \eqref{channel} $\e$-randomizing when $(U_i)$ are distributed according to any isotropic probability measure on $\mU(d)$ ?  We answer positively this question when $N \geq C d \log^6 d / \e^2$. This is our main result and appears as theorem \ref{th2}. Note that for non-Haar measures, previous results appearing in the literature \cite{hlsw,as,dn} involved the weaker trace-norm approximation $\| \Phi(\rho)-\rho_* \|_1 \leq \e$.

As opposed to the Haar measure, the measure $\mu$ need not have subgaussian tails, and we need more sophisticated tools to prove theorem \ref{th2}. We use recent results on suprema of empirical processes in Banach spaces. After early work by Rudelson \cite{rudelson} and Guédon--Rudelson \cite{gr}, a general sharp inequality was obtained by Guédon, Mendelson, Pajor and Tomczak--Jaegermann \cite{gmpt}. This inequality is valid in any Banach space with a sufficiently regular equivalent norm, such as $\ell_1^d$. The problem of $\e$-randomizing channels involves the supremum of an empirical process in the trace-class space $S_1^d$ (non-commutative analogue of $\ell_1^d$), which enters perfectly this setting.

The paper is organized as follows. Section \ref{section-background} contains background and precise statements of the theorems. Theorem \ref{th1} (for Haar measure) is proved in section \ref{section-th1}. Theorem \ref{th2} (for a general measure) is proved in section \ref{section-th2}. An appendix contains the needed facts about geometry and probability in Banach spaces.

{\bf Acknowledgement.} I thank Andreas Winter for several e-mail exchanges on the topic, and I am very grateful to Alain Pajor for showing me that the results of \cite{gmpt} can be applied here.

\section{Background and presentation of results}
\label{section-background}

Thoughout the paper, the letter $C$ and $c$ denote absolute constants whose value may change from occurrence to occurrence. We usually do not pay too much attention to the value of these constants.

\subsection{Schatten classes.} We write $\M(\C^d)$ for the space of complex $d \times d$ matrices. If $A \in \M(\C^d)$, let 
$s_1(A),\dots,s_d(A)$ denote the \emph{singular values} of $A$ (defined as the square roots of the eigenvalues of $AA^\dagger$). For $1 \leq p \leq \iy$, the \emph{Schatten $p$-norm} is defined as
\[ \|A\|_p = \left( \sum_{i=1}^d s_i(A)^p \right)^{1/p} .\]
For $p=\iy$, the definition should be understood as $\|A\|_{\iy} = \max s_i(A)$ and coincides with the usual operator norm. It is well-known (see \cite{bhatia}, section IV.2) that $(\M(\C^d),\|\cdot\|_p)$ is a complex normed space, denoted $S_p^d$ and called \emph{Schatten class}. The space $S_p^d$ is the non-commutative analogue of the space $\ell_p^d$. We write $B(S_p^d)$ for the unit ball of $S_p^d$.

The Schatten 2-norm (sometimes called Hilbert--Schmidt or Frobenius norm) is a Hilbert space norm associated to the inner product $\< A, B \> = \tr A^\dagger B$. This Hermitian structure allows to identify $\M(\C^d)$ with its dual space. Duality on Schatten norms holds as in the commutative case: if $p$ and $q$ are conjugate exponents (i.e. $1/p+1/q=1$), then the normed space dual to $S_p^d$ coincides with $S_q^d$.

\subsection{Completely positive maps}

We write $\Ms(\C^d)$ (resp. $\Mp(\C^d)$) for the set of self-adjoint (resp. positive semi-definite) $d \times d$ matrices. A linear map $\Phi : \M(\C^d) \to \M(\C^d)$ is said to \emph{preserve positivity} if $\Phi(\Mp(\C^d)) \subset \Mp(\C^d)$. Moreover, $\Phi$ is said to be \emph{completely positive} if for any $k \in \N$, the map
\[ \Phi \otimes \Id_{\M(\C^k)} : \M(\C^d \otimes \C^k) \to \M(\C^d \otimes \C^k) \]
preserves positivity. We use freely the canonical identification $\M(\C^d) \otimes \M(\C^k) \approx \M(\C^d \otimes \C^k)$.

If $(e_i)_{0 \leq i \leq d-1}$ denotes the canonical basis of $\C^d$, let $E_{ij} = \ket{e_i}\bra{e_j}$. To $\Phi : \M(\C^d) \to \M(\C^d)$ we associate $A_\Phi \in \M(\C^d \otimes \C^d)$ defined as
\[ A_\Phi = \sum_{i,j=1}^d E_{ij} \otimes \Phi(E_{ij}). \]
The matrix $A_\Phi$ is called the \emph{Choi matrix} of $\Phi$ ; it is well-known \cite{choi} that $\Phi$ is completely positive if and only if $A_\Phi$ is positive. Therefore, the set of completely positive operators on $\M(\C^d)$ is in one-to-one correspondence with $\Mp(\C^d \otimes \C^d)$. This correspondence is known as the \emph{Choi--Jamio\l kowski isomorphism}.

The spectral decomposition of $A_\Phi$ implies now the following: any completely positive map $\Phi$ on $\M(\C^d)$ can be decomposed as
\begin{equation} \label{kraus} \Phi : X \mapsto \sum_{i=1}^N V_iXV_i^\dagger. \end{equation}
Here $V_1,\dots,V_N$ are elements of $\M(\C^d)$. This decomposition is called a \emph{Kraus decomposition} of $\Phi$ of length $N$. The minimal length of a Kraus decomposition of $\Phi$ (called \emph{Kraus rank}) is equal to the rank of the Choi matrix $A_\Phi$. In particular it is always bounded by $d^2$.

\subsection{States and the completely depolarizing channel}

A \emph{state} on $\C^d$ is a element of $\Mp(\C^d)$ with trace 1. We write $\mD(\C^d)$ for the set of states ; it is a compact convex set with (real) dimension $d^2-1$. If $x \in \C^d$ is a unit vector, we write $P_x = \ket{x}\bra{x}$ for the associated rank one projector. The state $P_x$ is called a \emph{pure state}, and it follows from spectral decomposition that any state is a convex combination of pure states. A central role is played by the \emph{maximally mixed state} $\rho_* = \Id/d$ ($\rho_*$ is sometimes called the \emph{random state}).

A \emph{quantum channel} $\Phi : \M(\C^d) \to \M(\C^d)$ is a completely positive map which preserves trace: for any 
$X \in \M(\C^d),  \tr \Phi(X) = \tr X$.
Note that a quantum channel maps states to states. The trace-preserving condition reads on the Kraus decomposition \eqref{kraus} as
\[ \sum_{i=1}^N V_i^\dagger V_i = \Id. \]

An example of quantum channel that plays a central role in quantum information theory is the \emph{(completely) randomizing channel} (also called \emph{completely depolarizing channel}) $R : \M(\C^d) \to \M(\C^d)$.
\[ R : X \to \tr X \cdot \frac{\Id}{d} .\] 
The randomizing channel maps every state to $\rho_*$. The Choi matrix of $R$ is $A_R = \frac{1}{d} \Id_{\C^d \otimes \C^d}$. Since $A_R$ has full rank, any Kraus decomposition of $R$ must have length (at least) $d^2$. An explicit decomposition can be written using Fourier-type unitary operators: let $\omega = \exp(2i\pi/d)$ and $A$ and $B$ the matrices defined as
\begin{equation} \label{def-fourier} A(e_j) = e_{j+1 \textnormal{ mod }d} \ \ \ \ \ \ B(e_j) = \omega^j e_j .\end{equation}
For $1 \leq j,k \leq d$, define $V_{j,k}$ as the product $B^jA^k$. Note that $V_{j,k}$ belongs to the unitary group $\mU(d)$. A routine calculation (see also section \ref{section-isotropy}) shows that for any $X \in \M(\C^d)$,
\[ \frac{1}{d^2} \sum_{j,k=1}^d V_{j,k}X V_{j,k}^\dagger = \tr X \cdot \frac{\Id}{d} .\]
This is a Kraus decomposition of the randomizing channel.

\subsection{$\e$-randomizing channels}

We are interested in approximating the randomizing channel $R$ by channels with low Kraus rank. Following Hayden, Leung, Shor and Winter \cite{hlsw}, a quantum channel $\Phi$ is called \emph{$\e$-randomizing} if for any state $\rho \in \mD(\C^d)$, 
\[ \| \Phi(\rho)-\rho_* \|_{\iy} \leq \frac{\e}{d} .\]
It is equivalent to say that the spectrum of $\Phi(\rho)$ is contained in $[(1-\e)/d,(1+\e)/d]$ for any state $\rho$.
It has been proved in \cite{hlsw}  that there exist $\e$-randomizing channels with Kraus rank equal to $C d \log d/\e^2$ for some constant $d$. This is much smaller that $d^2$ (the Kraus rank of $R$). The construction is simple:  generate independent random Kraus operators according to the Haar measure on $\mU(d)$ and show that the induced quantum channel is $\e$-randomizing with nonzero probability. A key step in the proof is a discretization argument. We show that a simple trick improves the efficiency of the argument from \cite{hlsw} to prove the following

\begin{theo}[Haar-generated $\e$-randomizing channels]
\label{th1}
Let $(U_i)_{1 \leq i \leq N}$ be independent random matrices Haar-distributed on the unitary group $\mathcal{U}(d)$. Let $\Phi : \mathcal{\C}^d \to \mathcal{\C}^d$ be the quantum channel defined by
\[ \Phi(\rho) = \frac{1}{N} \sum_{i=1}^N U_i \rho U_i^\dagger .\]
Assume that $0<\e<1$ and $N \geq Cd/\e^2$. Then the channel $\Phi$ is $\e$-randomizing with nonzero probability.
\end{theo}

As often with random constructions, we actually prove that the conclusion holds true with \textit{large} probability: the probability of failure is exponentially small in $d$.

It is clear that the way $N$ depends on $d$ is optimal: if $\Phi$ is a $\e$-randomizing channel with $\e <1$, its Kraus rank must be at least $d$. This is because for any pure state $P_x$, $\Phi(P_x)$ must have full rank. 
The dependence in $\e$ is sharp for channels as constructed here, since lemma \ref{bernstein} below is sharp. However, it is not clear whether families of $\e$-randomizing channels with a better dependence in $\e$ can be found using a different construction, possibly partially deterministic.

One checks (using the value $c=1/6$ from \cite{hlsw} in lemma \ref{net} and optimizing over the net size) that the constant in theorem \ref{th1} can the chosen to, e.g.,  $C = 150$. This is presumably far from optimal.

\subsection{Isotropic measures on unitary matrices}
\label{section-isotropy}

Although the quantum channels constructed in theorem \ref{th1} have minimal Kraus rank, it can be argued that Haar-distributed random matrices are hard to generate in real-life situations. We introduce a wide class of measures on $\mU(d)$ that may replace the Haar measure. 

\begin{defn}
Say that a probability measure $\mu$ on $\mU(d)$ is \emph{isotropic} if for any $X \in \M(\C^d)$, 
\[ \int_{\mU(d)} UXU^\dagger d\mu(U) = \tr X \cdot \frac{\Id}{d} .\]
Similarly, a $\mU(d)$-valued random vector is called isotropic if its law is isotropic.
\end{defn}

\begin{lemma}
\label{lemma-isotropy}
Let $U=(U_{ij})$ be a $\mU(d)$-valued random vector. The following assertions are equivalent
\begin{enumerate}
 \item $U$ is isotropic.
 \item For any $X \in \M(\C^d)$, $\E | \tr UX^\dagger |^2 = \frac{1}{d} \|X\|_2^2$.
 \item For any indices $i,j,k,l$, $ \E U_{ij} \overline{U_{kl}} = \frac{1}{d} \delta_{i,k}\delta_{j,l}. $
\end{enumerate}
\end{lemma}

\begin{proof}
Implications $(3) \Rightarrow (1)$ and $(3) \Rightarrow (2)$ are easily checked by expansion. For $(1) \Rightarrow (3)$, simply take $X= \ket{e_j}\bra{e_k}$. Identity $(2)$ implies after polarization that for any $A,B \in \M(\C^d)$,
\[ \E \left[ \overline{\tr(UA^\dagger)} \tr(UB^\dagger) \right] = \frac{1}{d} \tr (AB^\dagger) ,\]
from which $(3)$ follows.
\end{proof}

Condition (3) of the lemma means that the covariance matrix of $U$ --- which is an element of $\M(\M(\C^d))$ --- is a multiple of the identity matrix. 

Of course the Haar measure is isotropic. Other examples are provided by discrete measures. Let $\mathscr{U} = \{ U_1,\dots,U_{d^2} \}$ be a family of unitary matrices, which are mutually orthogonal in the following sense: if $i \neq j$, then $\tr U_i^\dagger U_j=0$. For example, one can take $\mathscr{U} = \{B^jA^k\}_{1 \leq j,k \leq d}$, $A, B$ defined as \eqref{def-fourier}. Then the uniform probability measure on $\mathscr{U}$ is isotropic. Indeed, any $X \in \M(\C^d)$ can be decomposed as $X = \sum x_i U_i$ and condition (2) of lemma \ref{lemma-isotropy} is easily checked.

If we specialize to $d=2$, we obtain a random Pauli operator: assign probability $1/4$ to each of the following matrices
to get a isotropic measure
\[ \sigma_0 = \left( \begin{array}{cc} 1 & 0 \\ 0 & 1 \end{array} \right) 
, \ \ \ \sigma_1 = \left( \begin{array}{cc} 0 & 1 \\ 1 & 0 \end{array} \right) 
, \ \ \ \sigma_2 = \left( \begin{array}{cc} 0 & -i \\ i & 0 \end{array} \right) 
, \ \ \ \sigma_3 = \left( \begin{array}{cc} 1 & 0 \\ 0 & -1 \end{array} \right) .
\]

It is straightforward to check that isotropic vectors tensorize: if $X_1 \in \mU(d_1)$ and $X_2 \in \mU(d_2)$ are isotropic, so is $X_1 \otimes X_2 \in \mU(d_1d_2)$. If we work on $\M((\C^2)^{\otimes k})$, which corresponds to a set of $k$ qubits, a natural isotropic measure is therefore obtained by choosing independently a Pauli matrix on each qubit, i.e. assigning mass $1/4^k$ to the 
matrix $\sigma_{i_1} \otimes \cdots \otimes \sigma_{i_k}$ for any $i_1,\dots,i_k \in \{0,1,2,3\}^k$. 

\subsection{$\e$-randomizing channels for an isotropic measure}

We can now state our main theorem asserting that up to logarithmic terms, the Haar measure can be replaced in theorem 1 by simpler notions of randomness. We first state our result

\begin{theo}[General $\e$-randomizing channels]
\label{th2}
Let $\mu$ be an isotropic measure on the unitary group $\mU(d)$. Let $(U_i)_{1 \leq i \leq N}$ be independent $\mu$-distributed random matrices, and $\Phi : \mathcal{\C}^d \to \mathcal{\C}^d$ be the quantum channel defined as
\begin{equation} \label{defphi} \Phi(\rho) = \frac{1}{N} \sum_{i=1}^N U_i \rho U_i^\dagger .\end{equation}
Assume that $0<\e<1$ and $N \geq C d (\log d)^6/\e^2$. Then the channel $\Phi$ is $\e$-randomizing with nonzero probability.
\end{theo}

Theorem \ref{th2} applies in particular for product of random Pauli matrices as described in the previous section. It is of interest for certain cryptographic applications to know that $\e$-randomizing channels can be realized using Pauli matrices.

As opposed to theorem \ref{th1}, the conclusion of theorem \ref{th2} is not proved to hold with exponentially large probability. Applying the theorem with $\e \eta$ instead of $\e$ and using Markov inequality shows that $\Phi$ is $\e$-randomizing with probability larger than $1-\eta$ provided $N \geq C d \log^6 d / (\e^2\eta^2)$. 

Theorem \ref{th2} could be quickly deduced from a theorem appearing in \cite{gmpt}. However, the proof of \cite{gmpt} is rather intricate and uses Talagrand's majorizing measures in a central way. We give here a proof of our theorem which uses the simpler Dudley integral instead, giving the same result. We however rely an a entropy lemma from \cite{gmpt}, which appears as lemma \ref{entropy} in the appendix.

The $\log^6 d$ appearing in theorem \ref{th2} is certainly non optimal (see remarks at the end of the paper). However, some power of $\log d$ is needed, as shown by the next proposition.

\begin{pr}
Let $A,B$ defined as \eqref{def-fourier} and $\mu$ be the uniform measure on the set $\{B^jA^k\}_{1 \leq j,k \leq d}$. Consider $(X_i)$ independent $\mu$-distributed random unitary matrices. If the quantum channel $\Phi$ defined as \eqref{defphi} is $\frac{1}{2}$-randomizing with probability larger than $1/2$, then $N \geq c d \log d$.
\end{pr}

\begin{proof}
We will rely on the following standard result in elementary probability theory known as the coupon collector's problem (see \cite{durrett}, Chapter 1, example 5.10): if we choose independently and uniformly random elements among a set of $d$ elements, the mean (and also the median) number of choices before getting all elements at least once is equivalent to $d \log d$ for large $d$.

In our case, remember that $\omega = \exp(2i\pi/d)$ and for $0 \leq j \leq d-1$, define $x_j \in \C^d$ as
\[ x_j = \left( \frac{1}{\sqrt d}, \frac{\omega^j}{\sqrt d}, \frac{\omega^{2j}}{\sqrt d}, \dots , \frac{\omega^{(d-1)j}}{\sqrt d} \right)  .\]
Note that $\mB=(x_j)_{0 \leq j \leq d-1}$ is an orthonormal basis of $\C^d$ and that $B^jA^kx_0=x_j$. Consequently, if $U$ is $\mu$-distributed, the random state $UP_{x_0}U^\dagger$ equals $P_{x_j}$ with probability $1/d$. In the basis $\mB$, the matrix $\Phi(P_{x_0})$ is diagonal. Note that if $\Phi$ is $\frac{1}{2}$-randomizing, then $\Phi(P_{x_0})$ must have full rank. The reduction to the coupon collector's problem is now immediate.
\end{proof}

\section{Proof of theorem \ref{th1}: Haar-distributed unitary operators.}
\label{section-th1}

\noi 
The scheme of the proof is similar to \cite{hlsw}. We need two lemmas from there. The first is a deviation inequality sometimes known as Bernstein's inequality. The second is proved by a volumetric argument.

\begin{lemma}[Lemma II.3 in \cite{hlsw}]
\label{bernstein}
Let $\phi,\psi$ be pure states on $\C^d$ and $(U_i)_{1 \leq i \leq N}$ be independent Haar-distributed random unitary matrices. Then for every $0<\delta<1$,
\[ \P \left( \left| \frac{1}{N} \sum_{i=1}^N \tr (U_i\phi U_i^\dagger \psi) -\frac{1}{d} \right| \geq \frac{\delta}{d} \right) \leq 2 \exp (-c\delta^2 N)\]
\end{lemma}

\begin{lemma}[Lemma II.4 in \cite{hlsw}]
\label{net}
For $0<\delta<1$ there exists a set $\mathcal{N}$ of pure states on $\C^d$ with $|\mathcal{N}| \leq (5/\delta)^{2d}$, such that
for every pure state $\phi$ on $\C^d$, there exists $\phi_0 \in \mathcal{N}$ such that $\|\phi-\phi_0\|_1 \leq \delta$. Such a set $\mathcal{N}$ is called a $\delta$-net.
\end{lemma}

\noi The improvement on the result of \cite{hlsw} will follow from the next lemma

\begin{lemma}[Computing norms on nets]
 \label{norm-net} 
Let $\Delta : \mathcal{B}(\C^{d}) \to \mathcal{B}(\C^{d})$ be a Hermitian-preserving linear map. Let $A$ be the quantity 
\[ A = \sup_{\phi \in \mathcal{D}(\C^{d})} \| \Delta(\phi) \|_{\iy} = \sup_{\phi,\psi \in \mathcal{D}(\C^{d})} \left| \tr \psi \Delta(\phi) \right| \]
Let $0<\delta<1/2$ and $\mathcal{N}$ be a $\delta$-net as provided by lemma \ref{net}. We can evaluate $A$ as follows
\[ A \leq \frac{1}{1-2\delta} B ,\]
where
\[ B = \sup_{\phi_0,\psi_0 \in \mathcal{N}} \left| \tr \psi_0 \Delta(\phi_0) \right| \]
\end{lemma}

\begin{proof}[Proof of lemma \ref{norm-net}.]
First note that for any self-adjoint operators $a,b \in \mathcal{B}(\C^{d})$, we have
\begin{equation} \label{a1a2} \left| \tr b \Delta(a) \right| \leq A \|a\|_1 \|b\|_1. \end{equation}
By a convexity argument, the supremum in $A$ can be restricted to pure states. Given pure states $\phi,\psi \in \mathcal{D}(\C^{d})$, let $\phi_0,\psi_0 \in \mathcal{N}$ so that $\|\phi-\phi_0\|_1 \leq \delta$, $,\|\psi-\psi_0\|_1 \leq \delta$. Then
\[ \left| \tr \psi \Delta(\phi)\right| \leq \left| \tr (\psi-\psi_0) \Delta(\phi)\right| + \left| \tr \psi_0 \Delta(\phi-\phi_0)\right| + \left| \tr \psi_0 \Delta(\phi_0)\right|\] 
Using twice \eqref{a1a2} and taking supremum over $\phi,\psi$ gives $A \leq \delta A+\delta A+B$, hence the result.
\end{proof}

\begin{proof}[Proof of the theorem.] 
Let $R$ be the randomizing channel. Fix a $\frac{1}{4}$-net $\mathcal{N}$ with $|\mathcal{N}| \leq 20^{2d}$, as provided by lemma \ref{net}. Let $\Delta = R-\Phi$ and $A,B$ as in lemma \ref{norm-net}. 
Here $A$ and $B$ are random quantities and it follows from lemma \ref{norm-net} that
\[ \P \left(A \geq \frac{\e}{d} \right) \leq \P \left(B \geq \frac{\e}{2d} \right) .\]
Using the union bound and lemma \ref{bernstein}, we get
\[ \P \left(B \geq \frac{\e}{2d} \right) \leq 20^{4d} \cdot 2 \exp (-c\e^2 N/4).\]
This is less that $1$ if $N \geq Cd/\e^2$, for some constant $C$.
\end{proof}

\section{Proof of theorem \ref{th2}: general unitary operators.}
\label{section-th2}

A Bernoulli random variable is a random variable $\e$ so that $\P(\e=1)=\P(\e=-1)=1/2$. Recall that $C$ denotes an absolute constant whose value may change from occurrence to occurrence. We will derive theorem \ref{th2} from the following lemma.

\begin{lemma}
\label{rudelson}
Let $U_1,\dots,U_N \in \mU(d)$ be \emph{deterministic} unitary operators and let $(\e_i)$ be a sequence of independent Bernoulli random variables. Then
\begin{equation}
\label{eq-rudelson}
 \E_\e \sup_{\rho \in \mD(\C^d)} \left\| \sum_{i=1}^N \e_i U_i\rho U_i^\dagger \right\|_{\iy} \leq C (\log d)^{5/2} \sqrt{\log N} 
\sup_{\rho \in \mD(\C^d)} \left\| \sum_{i=1}^N U_i\rho U_i^\dagger \right\|_{\iy}^{1/2}. 
\end{equation}
\end{lemma}

\begin{proof}[Proof of theorem \ref{th2} (assuming lemma \ref{rudelson})]
Let $\mu$ be an isotropic measure on $\mU(d)$ and $(U_i)$ be independent $\mu$-distributed random unitary matrices.
Let $M$ be the random quantity
\[ M = \sup_{\rho \in \mD(\C^d)} \left\| \frac{1}{N} \sum_{i=1}^N U_i\rho U_i^\dagger - \frac{\Id}{d} \right\|_{\infty} \]
We are going to show that $\E M$ is small. 
The first step is a standard symmetrization argument. Let $(U'_i)$ be independent copies of $(U_i)$ and $(\e_i)$ be a sequence of 
independent Bernoulli random variables. We explicit as a subscript the random variables with repsect to which expectation is taken
\begin{eqnarray*}
\E M &\leq&  \E_{U,U'} \sup_{\rho \in \mD(\C^d)} \left\| \frac{1}{N} \sum_{i=1}^N U_i\rho U_i^\dagger - U_i'\rho U_i'^\dagger \right\|_{\iy} \\
 &=& \E_{U,U',\e} \sup_{\rho \in \mD(\C^d)} \left\| \frac{1}{N} \sum_{i=1}^N \e_i(U_i\rho U_i^\dagger - U_i'\rho U_i'^\dagger)
 \right\|_{\iy} \\
 & \leq & 2\E_{U,\e} \sup_{\rho \in \mD(\C^d)} \left\| \frac{1}{N} \sum_{i=1}^N \e_i U_i\rho U_i^\dagger \right\|_{\infty}
\end{eqnarray*}
The inequality of the first line is Jensen's inequality for $\E_{U'}$, while the equality on the second line holds since the distribution of $\rho \mapsto U_i\rho U_i^\dagger - U_i'\rho U_i'^\dagger$ is symmetric (as a $\M(\M(\C^d),\M(\C^d))$-valued random vector). We then decouple the expectations using lemma \ref{rudelson} for fixed $(U_i)$.
\begin{eqnarray*} \E M &\leq& \frac{C}{\sqrt{N}} (\log d)^{5/2} \sqrt{\log N} \E \sup_{\rho \in \mD(\C^d)} \left\| \frac{1}{N} \sum_{i=1}^N U_i\rho U_i^\dagger \right\|_{\iy}^{1/2} \\
& \leq & \frac{C}{\sqrt{N}} (\log d)^{5/2} \sqrt{\log N} \E \sqrt{M+\frac{1}{d}} \\
& \leq & \frac{C}{\sqrt{N}} (\log d)^{5/2} \sqrt{\log N} \sqrt{\E M+\frac{1}{d}} 
\end{eqnarray*}
Using the elementary implication
\[ X \leq \alpha \sqrt{X+\beta} \Longrightarrow   X \leq \alpha^2 +\alpha \sqrt{\beta}     \]
we find that $\E M \leq \e/d$ provided $N \geq C d \log^6 d/\e^2$.
\end{proof}

It remains to prove lemma \ref{rudelson}. We will use several standard concepts from geometry and probability in Banach spaces. All the relevant definitions and statements are postponed to the next section.

\begin{proof}[Proof of lemma \ref{rudelson}]
Let $Z$ be the quantity appearing in the left-hand side of \eqref{eq-rudelson}. By a convexity argument, the supremum  is attained for an extremal $\rho$, i.e. a pure state $P_x=\ket{x}\bra{x}$ for some unit vector $x$. Since the operator norm itself can be written as a supremum over unit vectors, we get
\[ Z
= \sup_{|x|=|y|=1} \left| \sum_{i=1}^N \e_i | \bra{y} U_i \ket{x} |^2 \right| 
= \sup_{|x|=|y|=1} \left| \sum_{i=1}^N \e_i | \tr U_i\ket{x}\bra{y}|^2 \right|
\leq \sup_{A\in B(S_1^d)} \left| \sum_{i=1}^N \e_i | \tr U_iA |^2 \right|.
\]
The last inequality follows from the fact that $B(S_1^d) = \conv \{ \ket{x}\bra{y}, |x|=|y|=1 \}$.
Let $\Phi : B(S_1^d) \to \R^N$ defined as 
\[ \Phi(A) = (|\tr U_1A|^2,\dots,|\tr U_NA|^2) .\]
We now apply Dudley's inequality (theorem \ref{dudley} in the next section) with $K=\Phi(B(S_1^d))$ to estimate $\E Z$ using covering numbers. This yields
\[ \E Z\leq C \int_0^{\iy} \sqrt{\log N(\Phi(B(S_1^d)),|\cdot|,\e)} d\e \]
where $|\cdot|$ denotes the Euclidean norm on $\R^N$.
Define a distance $\delta$ on $B(S_1^d)$ as
\[ \delta(A,B) = |\Phi(A)-\Phi(B)| = \left(\sum_{i=1}^N \left| |\tr U_iA|^2-|\tr U_iB|^2 \right|^2\right)^{1/2} .\]
We are led to the estimate 
\[ \E Z \leq C \int_0^{\iy} \sqrt{\log N(B(S_1^d),\delta,\e)} d\e.\]
Using the inequality $\left| |a|^2-|b|^2 \right| \leq |a-b|\cdot|a+b|$, the metric $\delta$ can be upped bounded as follows
\[ \delta(A,B)^2 \leq  \left( \sum_{i=1}^N |\tr U_i(A+B)|^2 \right) \sup_{1 \leq i \leq N} |\tr U_i(A-B)|^2 . \]
Let us introduce a new norm $|||\cdot|||$ on $\M(\C^d)$
\[ |||A||| = \sup_{1 \leq i \leq N} |\tr U_iA| .\]
Let $\theta$ be the number equal to 
\[ \theta := \sup_{A\in B(S_1^d)} \sum_{i=1}^N |\tr U_iA|^2 = \sup_{\rho \in \mD(\C^d)} \left\| \sum_{i=1}^N U_i\rho U_i^\dagger \right\|_{\iy}.\]
We get that for $A,B \in B(S_1^d)$, $\delta(A,B) \leq 2 \theta |||A-B|||$, and therefore
\[ \E Z \leq C\theta \int_0^{\iy} \sqrt{\log N(B(S_1^d),|||\cdot|||,\e)} d\e .\]
It remains to bound this new entropy integral. We split it into three parts, for $\e_0$ to be determined.
If $\e$ is large ($\e>1$), since $\|U_i\|_{\iy} = 1$, we get that $|||\cdot||| \leq \|\cdot\|_1$. This means that $N(B(S_1^d),|||\cdot|||,\e)=1$ and the integrand is zero.
If $\e$ is small ($0 < \e < \e_0$), we use the volumetric argument of lemma \ref{volumetric} 
\[ N(B(S_1^d),|||\cdot|||,\e) \leq N(B(S_1^d),\|\cdot\|_1,\e) \leq (3/\e)^{2d^2} .\]
In the intermediate range ($\e_0 \leq \e \leq 1$), let $q=\log d$ and $p = 1 + 1/(\log d-1)$ be the conjugate exponent. We are going to approximate the Schatten 1-norm by the Schatten $p$-norm. It is elementary to check that for $A \in \M(\C^d)$, $\|A\|_q \leq e \|A\|_\iy$. By dualizing 
\[ \|A\|_1 \leq e \|A\|_p \Longrightarrow N(B(S_1^d),|||\cdot|||,\e) \leq N(B(S_p^d),|||\cdot|||,\e/e).\]
We are now in position to apply lemma \ref{entropy} to the space $E = S_p^d$. By theorems \ref{th-2-convexity} and \ref{th-type2}, the 2-convexity constant of $S_p^d$ and the type 2 constant of $S_q^d$ (see next section for definitions) are bounded as follows
\[ T_2(S_q^d) \leq \lambda(S_p^d) \leq \sqrt{q-1} \leq \sqrt{\log d} .\]
Since $\|U_i\|_q \leq e$, the inequality given by lemma \ref{entropy} is
\[ \sqrt{\log N(B(S_1^d),|||\cdot|||,\e)} \leq \frac{C}{\e} (\log d)^{3/2} \sqrt{\log N}. \]
We now gather all the estimations
\[ \int_0^\iy \sqrt{\log N(B(S_1^d),|||\cdot|||,\e)} d\e\leq \int_0^{\e_0} \sqrt{2d^2} \log(3/\e) d\e+ C (\log d)^{3/2}\sqrt{\log N} \int_{\e_0}^1 \frac{1}{\e} d\e .\]
Choosing $\e_0=1/d$, an immediate computation shows that 
\[ \int_0^\iy \sqrt{\log N(B(S_1^d),|||\cdot|||,\e)} d\e\leq C (\log d)^{5/2} \sqrt{\log N}.\]
 This concludes the proof of the lemma.
\end{proof}

\section*{Appendix : Geometry of Banach spaces} \label{appendix}

In this last section, we gather several definitions and results from geometry and probability in Banach spaces.
We denote by $(E,\|\cdot\|)$ a real or complex Banach space (actually, in our applications $E$ will be finite-dimensional). We denote by $(E^*,\|\cdot\|_*)$ the dual Banach space.  

\subsection{Covering numbers}

\begin{Adefn}
If $(K,\delta)$ is a compact metric space, the \emph{covering number} or \emph{entropy number} $N(K,\delta,\e)$ is defined to be the smallest cardinality $M$ of a set $\{ x_1,\dots,x_M \} \subset K$ so that
\[ K \subset \bigcup_{i=1}^M B(x_i,\e) \]
where $B(x,\e) = \{ y \in K \st \delta(x,y) \leq \e \}$.
\end{Adefn}

An especially important case is when $K$ is a subset of $\R^n$ and $\delta$ is induced by a norm.
The next lemma is proved by a volumetric argument (see \cite{led-tal}, Lemma 9.5).

\begin{Alemma} \label{volumetric}
If $\|\cdot\|$ is a norm on $\R^n$ with unit ball $K$, then for every $\e>0$, $N(K,\|\cdot\|,\e) \leq (1+2/\e)^n$.
\end{Alemma}

The following theorem gives upper bounds on Bernoulli averages involving covering numbers. For a proof, see Lemma 4.5 and Theorem 11.17 in \cite{led-tal}.

\begin{Atheo}[Dudley's inequality]
\label{dudley}
Let $(\e_i)$ be independent Bernoulli random variables and $K$ be a compact subset of $\R^n$. Denote by $(x_1,\dots,x_n)$ the coordinates of a vector $x \in \R^n$. Then for some absolute constant $C$,
\[ \E \max_{x \in K} \sum_{i=1}^n \e_i x_i \leq 
C \int_0^\iy \sqrt{\log N(K,|\cdot|,\e)} d\e \]
where $|\cdot|$ denotes the Euclidean norm on $\R^n$.
\end{Atheo}

\subsection{2-convexity}

\begin{Adefn}
\label{def-2-convexity}
 A Banach space $(E,\|\cdot\|)$ is said to be \emph{2-convex} with constant $\lambda$ if for any $y,z 
\in E$, we have
\[ \| y \|^2 + \lambda^{-2} \|z\|^2 \leq \frac{1}{2} ( \left\|y+z \right\|^2 +\left\|y-z\right\|^2) .\]
The smallest such $\lambda$ is called the 2-convexity constant of $E$ and denoted by $\lambda(E)$.
\end{Adefn}

We say shortly that ``$E$ is 2-convex'' while the usual terminology should be ``$E$ has a modulus of convexity of power type 2''. This should not be confused with the notion of 2-convexity for Banach lattices \cite{lt}.

It follows from the parallelogram identity that a Hilbert space is 2-convex with constant 1. Other examples are $\ell_p$ and $S_p^d$ for $1 < p \leq 2$. The next theorem has been proved by Ball, Carlen and Lieb \cite{bcl}, refining on early work by Tomczak--Jaegermann \cite{tj}.

\begin{Atheo}
\label{th-2-convexity}
For $p \leq 2$, the following inequality holds for $A,B \in \M(\C^d)$
\[ \|A\|_p^2 + (p-1)\|B\|_p^2 \leq \frac{1}{2}\left( \|A+B\|_p^2 + \|A-B\|_p^2 \right).\]
Therefore, $S_p^d$ is 2-convex with constant $1/\sqrt{p-1}$.
\end{Atheo}

This property nicely dualizes. Indeed, it is easily checked (see \cite{bcl}, lemma 5) that $E$ is $2$-convex with constant $\lambda$ if any only if, for every $y,z \in E^*$, 
\[ \| y \|_*^2 + \lambda^{2} \|z\|_*^2 \geq \frac{1}{2} ( \left\|y+z \right\|_*^2 +\left\|y-z\right\|_*^2) .\]
In this case, $E^*$ is said to be \emph{2-smooth} with constant $\lambda$.

\subsection{Type 2}

\begin{Adefn}
 A Banach space $(E,\|\cdot\|)$ is said to have \emph{type 2} if there exists a constant $T_2$ so that for any finite sequence
 $y_1,\dots,y_N$ of vectors of $E$, we have
\begin{equation} \label{def-type2} \left( \E \left\| \sum_{i=1}^N \e_i y_i \right\|^2 \right)^{1/2} \leq T_2 
\left( \sum_{i=1}^N \|y_i\|^2 \right)^{1/2} .\end{equation}
The smallest possible $T_2$ is called the type 2 constant of $E$ and denoted $T_2(E)$.
Here, the expectation $\E$ is taken with respect to a sequence $(\e_i)$ of independent Bernoulli random variables.
\end{Adefn}

It follows from the (generalized) parallelogram identity that a Hilbert space has type 2 with constant $1$, and there is actually equality in \eqref{def-type2}. If a Banach space $E$ is 2-convex, then $E^*$ is 2-smooth. It is easily checked (by induction on the number of vectors involved) that a 2-smooth Banach space has type 2 with the same constant. We therefore have the inequality $T_2(E^*) \leq \lambda(E)$. In particular, theorem \ref{th-2-convexity} implies the following result, first proved by Tomczak-Jaegermann \cite{tj} with a worse constant.

\begin{Atheo}
\label{th-type2}
If $q \geq 2$, then $S_q^d$ has type 2 with the estimate
\[ T_2(S_q^d) \leq \sqrt{q-1} .\]
\end{Atheo}

\subsection{An entropy lemma}

The following lemma plays a key role in our proof. It appears as Lemma 1 in \cite{gmpt}. 

\begin{Alemma}
\label{entropy}
Let $E$ be a Banach space with unit ball $B(E)$. Assume that $E$ is 2-convex with constant $\lambda(E)$. 
Let $x_1,\dots,x_N$ be elements of $E^*$, and define a norm $|||\cdot|||$ on $E$ as
\[ |||y||| = \max_{1 \leq i \leq N} |  x_i(y) | .\]
Then for any $\e >0$ we have for some absolute constant $C$
\begin{equation} \label{eq-ent}
 \e \sqrt{\log N(B(E),|||\cdot|||,\e)} \leq C \lambda(E)^2 T_2(E^*) \sqrt{\log N} \max_{1\leq i \leq N} \|x_i\|_{E^*}.
\end{equation}
\end{Alemma}

The proof of lemma \ref{entropy} is based on a duality argument for covering numbers coming from \cite{bpst}. A positive answer to the duality conjecture for covering numbers (see \cite{amst} for a statement of the conjecture and recent results) would imply that the inequality \eqref{eq-ent} is valid without the factor $\lambda(E)^2$. This would improve our estimate in theorem \ref{th2} to $N \geq C d (\log d)^4/\e^2$.


\begin{thebibliography}{10}

\nocite{}

\bibitem{amtw}
A.~Ambainis, M.~Mosca, A.~Tapp and R.~de~Wolf, {\it Private quantum channels}. 41st Annual Symposium on Foundations of Computer Science (Redondo Beach, CA, 2000), 547--553, IEEE Comput. Soc. Press.

\bibitem{as}
A.~Ambainis and A.~Smith, {\it Small Pseudo-random Families of Matrices: Derandomizing Approximate Quantum Encryption}, Proceedings of RANDOM'04, 249--260.

\bibitem{amst}
S.~Artstein, V.~Milman, S.~Szarek and N.~Tomczak--Jaegermann, {\it On convexified packing and entropy duality}.  Geom. Funct. Anal.  {\bf 14}  (2004),  no. 5, 1134--1141.

\bibitem{bcl}
K.~Ball, E.~Carlen and E. Lieb, {\it Sharp uniform convexity and smoothness inequalities for trace norms}. 
Invent. Math. {\bf 115} (1994), no. 3, 463--482. 

\bibitem{bhatia}
R. Bhatia, {\it Matrix analysis}. Graduate Texts in Mathematics {\bf 169}. Springer-Verlag, 1997.

\bibitem{bpst}
J.~Bourgain, A.~Pajor, S.~Szarek and N.~Tomczak-Jaegermann, {\it On the duality problem for entropy numbers of operators}.  Geometric aspects of functional analysis (1987--88),  50--63, Lecture Notes in Math. {\bf 1376} (1989).

\bibitem{choi}
M.~D.~Choi, {\it Completely positive linear maps on complex matrices}. Linear Algebra and Appl. {\bf 10}  (1975), 285--290.

\bibitem{dn}
P.~Dickinson and A.~Nayak, {\it Approximate Randomization of Quantum States With Fewer Bits of Key},
AIP Conference Proceedings {\bf 864}, 18--36 (2006).

\bibitem{durrett}
R.~Durrett, {\it Probability. Theory and examples}, The Wadsworth \& Brooks/Cole Statistics/Probability Series, 1991.

\bibitem{gmpt} 
O.~Guédon, S.~Mendelson, A.~Pajor and N.~Tomczak--Jaegermann, {\it Majorizing measures and proportional subsets of bounded orthonormal systems}, preprint (2008).

\bibitem{gr}
O.~Guédon and M.~Rudelson,
{\it $L_p$-moments of random vectors via majorizing measures}, Adv. Math.  {\bf 208}  (2007),  no. 2, 798--823.

\bibitem{hlsw}
P.~Hayden, D.~Leung, P.~W.~Shor and A.~Winter, {\it Randomizing quantum states: constructions and applications}, Comm. Math. Phys. \textbf{250} (2004), 371--391. 

\bibitem{led-tal}
M.~Ledoux and M.~Talagrand, {\it Probability in Banach spaces. Isoperimetry and processes}. Ergebnisse der Mathematik und ihrer Grenzgebiete (3) {\bf 23}. Springer-Verlag, 1991.

\bibitem{lt}
J.~Lindenstrauss and L.~Tzafriri, {\it Classical Banach spaces. II. Function spaces}. Ergebnisse der Mathematik und ihrer Grenzgebiete {\bf 97}. Springer-Verlag, 1979.

\bibitem{rudelson}
M.~Rudelson, {\it Random vectors in the isotropic position}. J. Funct. Anal. {\bf 164}  (1999),  no. 1, 60--72.

\bibitem{tj}
N.~Tomczak-Jaegermann, {\it The moduli of smoothness and convexity and the Rademacher averages of trace classes $S_p(1\leq p<\infty )$.} Studia Math. {\bf 50}  (1974), 163--182.

\bibitem{winter}
A.~Winter, {\it The maximum output $p$-norm of quantum channels is not multiplicative for any $p$> 2}, arxiv \verb!0707.0402!

\end{thebibliography}
\end{document}